\RequirePackage{fix-cm}
\documentclass{svjour3}                     % onecolumn (standard format)
\def\makeheadbox{\relax}

\usepackage{natbib}
\usepackage{caption,subcaption}
\usepackage{graphicx}
\usepackage[utf8]{inputenc}
\usepackage{amsmath,amssymb,amsfonts}
\usepackage[title]{appendix}

\newcommand{\kone}{\langle k \rangle}
\newcommand{\ktwo}{\langle k^2 \rangle}
\newcommand{\kthree}{\langle k^3 \rangle}

\begin{document}

\title{A Low-Dimensional Network Model for an SIS Epidemic}
% Grants or other notes about the article that should go on the front
% page should be placed within the \thanks{} command in the title
% (and the %-sign in front of \thanks{} should be deleted)
%
% General acknowledgments should be placed at the end of the article.

\subtitle{Analysis of the Super Compact Pairwise Model}

%\titlerunning{Short form of title}        % if too long for running head

\author{Carl Corcoran         \and
        Alan Hastings
}

%\authorrunning{Short form of author list} % if too long for running head

\institute{C. Corcoran \at
              Department of Mathematics, University of California, Davis, USA \\
              \email{ctcorcoran@ucdavis.edu}           %  \\
%             \emph{Present address:} of F. Author  %  if needed
           \and
           A. Hastings \at
              Department of Environmental Science and Policy, University of California, Davis, USA and Santa Fe Institute, Santa Fe, USA
}

\date{\today}
% The correct dates will be entered by the editor

\maketitle

\begin{abstract}
Network-based models of epidemic spread have become increasingly popular in recent decades. Despite a rich foundation of such models, few low-dimensional systems for modeling SIS-type diseases have been proposed that manage to capture the complex dynamics induced by the network structure. We analyze one recently introduced model and derive important epidemiological quantities for the system. We derive the epidemic threshold and analyze the bifurcation that occurs, and we use asymptotic techniques to derive an approximation for the endemic equilibrium when it exists. We consider the sensitivity of this approximation to network parameters, and the implications for disease control measures are found to be in line with the results of existing studies. 

% \PACS{PACS code1 \and PACS code2 \and more}
% \subclass{MSC code1 \and MSC code2 \and more}
\end{abstract}

%%%%%%%%%%%%%%%%%%%%%%%%%%%%%%%%%%%%%%%%%%%%%%%%%%%%%%%%%%%%%%%%%%%%%%%%%%%%%%%
%%%%%%%%%%%%%%%%%%%%%%%%%%%%%%%%%%%%%%%%%%%%%%%%%%%%%%%%%%%%%%%%%%%%%%%%%%%%%%%
\section{Introduction} \label{intro}

In the past few decades, network-based models of epidemic spread have become a central topic \citep{kiss_mathematics_2017,pastor-satorras_epidemic_2015} in epidemiology. Their ability to capture mathematically the complex structure of transmission interactions makes them an invaluable theoretical paradigm. Mathematically, a network is modeled as a graph consisting of a set of nodes that are connected by a set of links (called edges). In the context of epidemiology, typically nodes represent individuals, and edges represent interactions that can transmit the infection. Used in conjunction with compartment models, the disease natural history determines the number of possible states an individual node might be in at any point in time. When disease spread is modeled as a continuous time Markov chain, the network size and disease natural history can lead to high dimensional state spaces. For example, in a network with $N$ nodes where individual nodes can be in $m$ possible states, the size of the state space for the network is $m^N.$ Efforts to describe this process with a system of ordinary differential equations are similarly hampered by size---the Kolmogorov equations governing this system are exact, but prohibitively large. Thus, an important goal in network-based modeling has been to find (relatively) low-dimensional system that accurately approximates the underlying high-dimensional system.

Many approaches \citep{pastor-satorras_epidemic_2001-1,pastor-satorras_epidemic_2015,miller_edge-based_2012,barnard_epidemic_2019,karrer_message_2010} in recent years have sought to introduce models with systems of a manageable size. Pairwise models \citep{keeling_effects_1999-1,eames_modeling_2002,house_insights_2011} have been a popular and fruitful approach to this question. Derived from the Kolmogorov equations and exact in their unclosed form \citep{taylor_markovian_2012}, pairwise models consider the evolution of not just the expected number of nodes in a given state, but also pairs and triples of nodes. The dynamical variables are of the form $[A]$ (the expected number of nodes in state $A$), $[AB]$ (the expected number of pairs in state $A-B$), and $[ABC]$ (the expected number of triples in state $A-B-C$). Higher-order groupings of nodes are also considered but rarely written, as dimension-reduction efforts often focus on approximating the expected number of triples in terms of pairs and individuals nodes. 
Pairwise models have been successful with a variety of different network types, with models developed for networks with heterogeneous degree \citep{eames_modeling_2002}, weighted networks \citep{rattana_class_2013}, directed networks \citep{sharkey_pair-level_2006}, and networks with motifs \citep{house_motif-based_2009,keeling_systematic_2016} to name a few. Moreover, pairwise models have been developed for a variety of disease natural histories, with particular focus on SIR (susceptible-infected-recovered) and SIS (susceptible-infected-susceptible) models.

In this paper, we consider an SIS pairwise model for networks with heterogeneous degree. SIS dynamics are used to model diseases where no long term immunity is conferred upon recovery, leading to their frequent application to sexually transmitted infections such as chlamydia or gonorrhea \citep{eames_modeling_2002}. 
Contact networks for diseases of this type frequently involve heterogeneity in the number of contacts for individuals, and thus node degree becomes an essential concept. The degree of a node in a network is the number of edges to which the node is connected, and thus the number of potential infectious contacts. In this way, heterogeneous networks can capture complex disease dynamics. An essential tool when working with such networks is the degree distribution, defined by $p_k$ which is the probability a randomly selected node has degree $k.$ The degree distribution has played an important role in dimension reduction approximations for pairwise models.

For the SIR-type diseases, accurate low-dimensional models have been derived from the pairwise family using probability generating functions \citep{miller_edge-based_2012}, complete with conditions for finding the final size of the epidemic. Despite the successes of the SIR case, the dimension reduction techniques in \citet{miller_edge-based_2012} do not apply to the SIS case. Instead, the development of low-dimensional models of SIS-type disease spread on networks have relied on moment closure approximations. Under the assumption of a heterogeneous network with no clustering, \citet{house_insights_2011} introduced an approximation reducing the system size from $\mathcal{O}(N^2)$ to $\mathcal{O}(N)$, where $N$ is the number of nodes in the network. Termed the compact pairwise model (CPW), it has shown good agreement with stochastic simulations despite its considerably smaller size. However, the number of model equations still grows as the maximum degree of the network, making its application challenging for large networks with significant degree heterogeneity. 
Perhaps the most successful model in reducing the number of equations of the CPW for SIS-type diseases is the super compact pairwise model (SCPW) \citep{simon_super_2016-1}. The system consists of only four equations, with network structure being encoded to the model through the first three moments of the degree distribution. While Simon and Kiss demonstrated excellent agreement between the CPW and the SCPW, bifurcation analysis of the model and an explicit formula for the endemic steady state remain to be done.

This paper sets out on that analysis of the SCPW model. A common point of investigation among models of SIS-type diseases is the disease-free equilibrium (DFE) that loses stability as a relevant parameter passes a critical value known as the epidemic threshold \citep{pastor-satorras_epidemic_2001-1,pastor-satorras_epidemic_2002,boguna_epidemic_2002}. The epidemic threshold serves as a dividing point between two qualitatively different types of outbreaks. Below the epidemic threshold, any outbreak will die out; above the epidemic threshold, the system converges asymptotically to a stable equilibrium where the disease remains endemic in the population. Many studies follow the ``next generation matrix" approach for the basic reproduction number $R_0$ \citep{van_den_driessche_reproduction_2002}  to characterize the epidemic threshold. We follow a more conventional bifurcation analysis to derive the epidemic threshold and offer a proof that the system undergoes a transcritical bifurcation, as one might expect. Perhaps more importantly, the SCPW's small fixed number of equations presents an excellent opportunity to investigate the endemic equilibrium for SIS models on heterogeneous networks, which has been heretofore inhibited by large system size. We present a novel asymptotic approach to approximating the endemic equilibrium, leveraging the low-dimensionality of the model. The results presented further our understanding of the SCPW model specifically, and suggest potential new avenues in the challenging problem of analytically determining the nontrivial steady state of pairwise models of SIS-type diseases.

The paper is structured as follows: in Section \ref{model}, we nondimensionalize the model and reduce the number of equations to 3 to facilitate computations. In Section \ref{threshold}, we derive the epidemic threshold and show that the system undergoes a forward transcritical bifurcation. In Section \ref{asymptotics}, we tackle the endemic steady state that emerges through the bifurcation. We use asymptotic methods to approximate the size of the endemic steady state under two regimes---the system near the epidemic threshold and the system far away from the epidemic threshold---and give examples of the efficacy of these approximations on prototypical networks. Finally, we examine the implications of these two approximations. In line with existing studies \citep{eames_modeling_2002}, we find that control measures for reducing the prevalence at the endemic equilibrium may require different tactics depending on the regime.    
%%%%%%%%%%%%%%%%%%%%%%%%%%%%%%%%%%%%%%%%%%%%%%%%%%%%%%%%%%%%%%%%%%%%%%%%%%%%%%%%%%
%%%%%%%%%%%%%%%%%%%%%%%%%%%%%%%%%%%%%%%%%%%%%%%%%%%%%%%%%%%%%%%%%%%%%%%%%%%%%%%%%%
\section{Model} \label{model}

The essential characteristic of pairwise models of SIS epidemics is dynamical equations for not just the expected number of nodes in each state, but also pairs and triples of nodes. At the node level, $[S]$ and $[I]$ are the expected number of susceptible and infectious nodes respectively. At the pair level, $[SI]$ is the expected number of connected pairs of susceptible and infectious nodes, while $[SS]$ and $[II]$ are the expected numbers of connected susceptible-susceptible and infectious-infectious pairs respectively. The full pairwise model further requires equations for the expected number of triples ($[SSI]$ and $[ISI]$) and higher motifs as well:
\begin{align*}
    \dot{[S]} &= \gamma [I] - \tau [SI], \\
    \dot{[I]} &= \tau [SI] - \gamma [I], \\
    \dot{[SI]} &= \gamma([II] - [SI]) +\tau ([SSI]-[ISI]-[SI]), \\
    \dot{[SS]} &= 2\gamma[SI]-2\tau [SSI], \\
    \dot{[II]} &= -2\gamma[II]+2\tau([ISI]+ [SI]).
\end{align*}
The CPW closes the system by approximating the expected number of triples as
\begin{equation*}
    [ASI] \approx [AS][SI]\frac{S_2-S_1}{S_1^2},
\end{equation*}
where $S_1$ and $S_2$ are the first and second moments of the distribution of susceptible nodes; that is
\begin{equation*}
    S_1 = \sum_k k[S_k] =[SS] + [SI],\:\: S_2 = \sum_k k^2[S_k],
\end{equation*}
where $[S_k]$ is the expected number of susceptible nodes with degree $k.$ Unfortunately $S_2$ cannot be expressed exactly in terms of $[S],[I],[SI],[SS],$ and $[II]$ only, so the SCPW model offers an approximation that depends on these variables and moments of the degree distribution.

The SCPW model derived in \citet{simon_super_2016-1} is given as
\begin{align}
    \dot{[S]} &= \gamma [I] - \tau [SI], \label{eq:1}\\
    \dot{[I]} &= \tau [SI] - \gamma [I] \label{eq:2}\\
    \dot{[SI]} &= \gamma([II] - [SI]) - \tau [SI]+ \tau [SI] ([SS]-[SI])Q, \label{eq:3}\\
    \dot{[SS]} &= 2\gamma[SI]-2\tau [SI][SS]Q, \label{eq:4}\\
    \dot{[II]} &= -2\gamma[II]+2\tau [SI]+2\tau [SI]^2Q,  \label{eq:5}
\end{align}
    where 
\begin{equation*}
Q = \frac{1}{n_S [S]}\left(\frac{\ktwo(\ktwo - \kone n_S)+\kthree (n_S-\kone)}{n_S(\ktwo - \kone^2)}-1\right),\: n_S = \frac{[SI]+[SS]}{[S]},
\end{equation*}
$\langle k^n \rangle$ is the $n$th moment of the degree distribution, $\tau$ is the transmission rate, and $\gamma$ is the recovery rate. Here, the quantity $Q$ serves as an approximation of $(S_2-S_1)/S_1^2.$ As well, the quantities $[S],[I],[SI],[SS],[II]$ satisfy  conservation equations
\begin{align}
[S] + [I] &= N, \label{eq:6}\\
2[SI] + [SS]+[II] &= \langle  k\rangle N .\label{eq:7}
\end{align}

With the goal of performing bifurcation and asymptotic analyses in mind, nondimensionalizing the SCPW model is a natural first step. To do so, we will rearrange the equations (\ref{eq:3})-(\ref{eq:5}) so that the network parameters $\kone, \ktwo, \kthree$ are consolidated into more workable constants. First, we rewrite $Q$ as
\begin{equation}
    Q = \frac{\alpha [S]}{([SI]+[SS])^2}+\frac{\beta}{[SI]+[SS]}, \label{eq:8}
\end{equation}
where 
\begin{equation}
    \alpha = \frac{\ktwo^2-\kone\kthree}{\ktwo - \kone^2},\:\: \beta = \frac{\kthree-\ktwo\kone}{\ktwo - \kone^2}-1. \label{eq:9}
\end{equation}
A natural nondimensionalization of this system is to scale the number of nodes and links in each state to the proportion of nodes and pairs in each state: $v = [S]/N, w = [I]/N, x=[SI]/(\kone N),y=[SS]/(\kone N),z=[II]/(\kone N).$ As well, a natural rescaling of time is $T = t/\gamma,$ which prompts the defining of the transmission-recovery rate ratio $\delta = \tau/\gamma.$ The introduction of $\delta$ consolidates the two epidemiological parameters $\tau$ and $\gamma$ into a single nondimensional parameter, so any changes to epidemiology of the disease will be captured in $\delta$ alone. With these substitutions, the system (\ref{eq:1})-(\ref{eq:5}) becomes
\begin{align}
    \dot{v} &= w - \kone\delta x, \label{eq:10}\\
    \dot{w} &= \kone \delta x - w, \label{eq:11}\\
    \dot{x} &= z - \left(\delta + 1\right)x+\frac{\alpha \delta}{\kone }\cdot\frac{vx(y-x)}{(x+y)^2}+\beta\delta\cdot \frac{x(y-x)}{x+y},\label{eq:12}\\
    \dot{y} &= 2x - \frac{2\alpha\delta}{\kone }\cdot\frac{vxy}{(x+y)^2}-2\beta\delta\cdot\frac{xy}{x+y},\label{eq:13}\\
    \dot{z} &= -2z+2\delta x+\frac{2\alpha\delta}{\kone}\cdot\frac{vx^2}{(x+y)^2}+2\beta\delta\cdot\frac{x^2}{x+y},\label{eq:14}
\end{align}
where the dot notation represents the derivative with respect to the nondimensional time variable $\frac{d}{dT}$. The conservation equations (\ref{eq:6}) and (\ref{eq:7}) become 
\begin{align}
v+w &= 1, \label{eq:15}\\
2x+y+z &= 1, \label{eq:16}
\end{align}
respectively.

At this point, the conservation equations can be used to reduce the system to a mere 3 equations. However, the elimination of different equations for different analyses will be convenient. For characterizing the bifurcation undergone by the disease-free equilibrium (DFE), it is  convenient to work with variables that are 0 at the DFE. For approximating the endemic steady state using asymptotic methods, the most parsimonious equations will make the algebraic manipulation required easier. Thus, we will work with slightly different (but equivalent) characterizations of (\ref{eq:10})-(\ref{eq:14}) in the sections that follow.

%%%%%%%%%%%%%%%%%%%%%%%%%%%%%%%%%%%%%%%%%%%%%%%%%%%%%%%%%%%%%%%%%%%%%%%%%%%%%%%%%%
%%%%%%%%%%%%%%%%%%%%%%%%%%%%%%%%%%%%%%%%%%%%%%%%%%%%%%%%%%%%%%%%%%%%%%%%%%%%%%%%%
\section{Epidemic Threshold} \label{threshold}

To derive the epidemic threshold, we consider the stability of the DFE in terms of the epidemiological parameter $\delta.$ We will show that as $\delta$ increases through a critical value $\delta_c,$ the DFE loses stability. Typically as the DFE loses stability,  an asymptotically stable endemic equilibrium emerges. The SCPW is no exception, and here we derive the epidemic threshold, with a proof that the system undergoes a transcritical bifurcation (and thus an endemic equilibrium emerges) when $\delta = \delta_c$ included in Appendix A.

First, we use the conservation equations (\ref{eq:15}) and (\ref{eq:16}) to eliminate equations (\ref{eq:10}) and (\ref{eq:13}). The resulting system is 
\begin{align}
    \dot{w} &= \kone \delta x - w, \label{eq:17}\\
    \dot{x} &= z - \left(\delta + 1\right)x+\frac{\alpha \delta}{\kone }\cdot\frac{(1-w)x(1-3x-z)}{(1-x-z)^2}+\beta\delta\cdot \frac{x(1-3x-z)}{1-x-z},\label{eq:18}\\
    \dot{z} &= -2z+2\delta x+\frac{2\alpha\delta}{\kone}\cdot\frac{(1-w)x^2}{(1-x-z)^2}+2\beta\delta\cdot\frac{x^2}{1-x-z}.\label{eq:19}
\end{align}
Though ostensibly a messier choice of equation reduction, we note that at the DFE, $[I] = [SI] = [II] =0,$ so $w=x=z=0.$ The notation
\begin{equation}
    \dot{\mathbf{x}} = \begin{bmatrix} \dot{w}\\ \dot{x} \\ \dot{z}\\ \end{bmatrix} = \begin{bmatrix} F_1(w,x,z)\\ F_2(w,x,z) \\ F_3(w,x,z) \\ \end{bmatrix} = \mathbf{F}(\mathbf{x}) \label{eq:20}
\end{equation}
will be convenient moving forward. To determine the stability of the DFE, we compute the Jacobian at $\mathbf{x} = \Vec{0}:$
\begin{equation}
    D\mathbf{F} = \begin{bmatrix}-1 & \kone \delta & 0 \\ 0 & \left(\dfrac{\alpha}{\kone}+\beta\right)\delta - (\delta+1)& 1 \\ 0 & 2\delta &-2\\  \end{bmatrix}.\label{eq:21}
\end{equation}
A straightforward computation shows that
\begin{equation}
    \frac{\alpha}{\kone} + \beta = \frac{\ktwo-\kone}{\kone} = \bar{k}. \label{eq:22}
\end{equation}
We can write $D\mathbf{F}$ as a block triangular matrix as
\begin{equation*}
    D\mathbf{F} = \begin{bmatrix}-1 & A \\ 0 & B\\ \end{bmatrix},
\end{equation*}
where the dimensions $A$ and $B$ respectively are $1\times 2$ and $2\times 2$. The properties of determinants of block matrices tell us that the eigenvalues of $D\mathbf{F}$ are $-1$ and the eigenvalues of $B,$ which will determine the stability of the DFE. 

We appeal here to the trace-determinant theorem, which tells us the eigenvalues $\xi$ of the $2\times 2$ matrix $B$ are given by
\begin{equation*}
    \xi = \frac{\text{Tr}(B)}{2} \pm \frac{\sqrt{(\text{Tr}(B))^2-4\,\text{Det}(B)}}{2}.
\end{equation*}
First, we observe that these eigenvalues are real, as
\begin{equation}
    \text{Tr}(B)^2-4\,\text{Det}(B) = (\delta(\bar{k}-1)+1)^2+8\delta, \label{eq:23} 
\end{equation}
which is clearly positive. As a consequence, for the DFE to be stable we must have $\text{Tr}(B) <0$ and $\text{Det}(B) >0.$ The determinant can be written 
\begin{equation}
    \text{Det}(B) = 2(1-\delta\bar{k}), \label{eq:24}
\end{equation}
and is thus positive if and only if $\delta < 1/\bar{k}.$ Moreover, if $\delta< 1/\bar{k},$ then 
\begin{equation*}
\text{Tr}(B) < (\bar{k}-1)/\bar{k}-3 = -2-1/\bar{k} <0.
\end{equation*}
Therefore, we conclude that the DFE is stable for $\delta< 1/\bar{k}$ and unstable for $\delta > 1/\bar{k}.$ Thus, the epidemic threshold is the critical value of the bifurcation parameter $\delta:$
\begin{equation}
    \delta_c = \frac{\kone}{\ktwo - \kone}.\label{eq:25}
\end{equation}
Notably, this threshold value is identical to that of the CPW as shown in \citet{kiss_mathematics_2017}. However it remains to be shown that a bifurcation actually does occur here, and that a asymptotically stable endemic steady state emerges. To prove this, we apply a theorem of \citet{castillo-chavez_dynamical_2004} in Appendix A.

%%%%%%%%%%%%%%%%%%%%%%%%%%%%%%%%%%%%%%%%%%%%%%%%%%%%%%%%%%%%%%%%%%%%%%%%%%%%%%%%%%
%%%%%%%%%%%%%%%%%%%%%%%%%%%%%%%%%%%%%%%%%%%%%%%%%%%%%%%%%%%%%%%%%%%%%%%%%%%%%%%%%%
\section{The Endemic Equilibrium} \label{asymptotics}

With existence of an endemic steady state established, we turn to the question of finding an approximate analytic expression for it. In general, this is a difficult proposition with epidemic models on networks owing to the frequently high-dimensional nature of the dynamical systems. An exact closed-form expression for the endemic equilibrium of the SCPW model requires solving a system of polynomial equations in multiple variables, which we do not attempt here. However, with asymptotic techniques, a workable approximation can be derived for two cases of $\delta$: near the epidemic threshold ($\delta \approx \delta_c$), and far away from it ($\delta >> \delta_c$). Two challenges are apparent. First, how to eliminate equations to facilitate asymptotic expansions of the equilibrium and second, the choice of small nondimensional parameter be in each case.

Unlike in Section \ref{threshold}, the most parsimonious characterization of (\ref{eq:10})-(\ref{eq:14}) is desirable. So we eliminate (\ref{eq:11}) and (\ref{eq:14}) with the conservation equations. To promote the finding of a small nondimensional parameter, we rewrite the resulting system using $\delta = \delta_c\cdot \frac{\delta}{\delta_c}$ and incorporate the constants $\sigma = \kone \delta_c, \lambda = \alpha \delta_c/\kone, \mu = \beta \delta_c.$ With these alterations, the system becomes
\begin{align}
    \dot{v} &= 1-v-\sigma\frac{\delta}{\delta_c}x, \label{eq:26} \\ 
    \dot{x} &= 1-y-\left(3+\delta_c\frac{\delta}{\delta_c}\right)x+\lambda\frac{\delta}{\delta_c}\frac{vx(y-x)}{(x+y)^2}+\mu\frac{\delta}{\delta_c}\frac{x(y-x)}{x+y}, \label{eq:27}\\ 
    \dot{y}&= 2x-2\lambda\frac{\delta}{\delta_c}\frac{vxy}{(x+y)^2}-2\mu\frac{\delta}{\delta_c}\frac{xy}{x+y}. \label{eq:28} 
\end{align}
At the endemic equilibrium, $\dot{v} = \dot{x} = \dot{y} =0.$ We can solve (\ref{eq:26}) for $v$ and substitute into (\ref{eq:27}) and (\ref{eq:28}). With some rearrangement of terms (and adding (\ref{eq:28}) to (\ref{eq:27})) we arrive at the system of polynomial equations that determines the endemic steady state:
\begin{align}
    0 &= \left(\frac{\delta_c}{\delta}\right)^2(1-y-2x)(x+y)^2+\frac{\delta_c}{\delta}\left(\delta_c x (x+y)^2+\lambda x^2 + \mu x(x+y)\right) \nonumber \\
    &\qquad +\lambda\sigma x^3 = P(x,y), \label{eq:29}\\ 0 &=\left(\frac{\delta_c}{\delta}\right)^2(x+y)^2 -\frac{\delta_c}{\delta}\left(\lambda y+\mu y (x+y)\right) + \lambda\sigma x y = Q(x,y)\label{eq:30}.
\end{align}
For the endemic steady state, we are interested in knowing the prevalence when the system is at equilibrium: $w^*$. We use the following procedure to approximate the solution.
\begin{enumerate}{}
    \item Express $\delta_c/\delta$ in terms of a small parameter.
    \item Use the Implicit Function Theorem to linearize $P(x,y)=0$ as 
    \begin{equation*}
        y \approx  \Tilde{y} - \dfrac{P_x(\Tilde{x},\Tilde{y})}{P_y(\Tilde{x},\Tilde{y})}(x-\Tilde{x})
    \end{equation*}
    around a point $(\Tilde{x},\Tilde{y})$ that is mathematically and/or biologically justified for the given regime.
    \item Expand $x, y$, and other relevant quantities in terms of the small parameter.
    \item Substitute the expansions into $Q(x,y)=0$ and obtain a regular perturbation problem and find an asymptotic solution for the equilibrium value $x$, which approximates $x^*$.
    \item Apply the relation $w^* = (\delta_c/\delta)^{-1}\sigma x^*$ to obtain an asymptotic series for the prevalence at the endemic equilibrium.
\end{enumerate}
We describe the results of this procedure for each case in the remainder of this section--the details of the computations are included in Appendix B.

%%%%%%%%%%%%%%%%%%%%%%%%%%%%%%%%%%%%%%%%%%%%%%%%%%%%%%%%%%%%%%%%%%%%%%%%%%%%%%%%%%
\subsection{Case 1: Near the epidemic threshold ($\delta \approx \delta_c$)}
For $\delta \approx \delta_c,$ we choose $\eta = 1-\delta_c/\delta$ as a small parameter. In terms of this small parameter, (\ref{eq:29}) and (\ref{eq:30}) become:
\begin{align}
    0 &= (1-\eta)^2(1-y-2x)(x+y)^2 \nonumber\\
    &\qquad -(1-\eta)\left(\delta_cx(x+y)^2+\lambda x^2 +\mu x^2(x+y)\right)+\lambda\sigma x^3, \label{eq:31}\\
    0 &= (1-\eta)^2(x+y)^2 -(1-\eta)\left(\lambda y +\mu y (x+y)\right) +\lambda\sigma xy. \label{eq:32}
\end{align}
When $\delta \approx \delta_c,$ an endemic steady state has just emerged, so we can view this equilibrium as a small perturbation to the steady state $x=0,\:y=1.$ 
Linearizing $P(x,y) = 0$ about this point gives 
    \begin{equation}
        y \approx 1 - \left(2+\frac{\delta_c}{1-\eta}\right)x.\label{eq:33}
    \end{equation}
Expanding 
    \begin{align}
        2+\frac{\delta_c}{1-\eta} = 2+\delta_c(1+\eta+\eta^2+\mathcal{O}(\eta^3)), \label{eq:34}\\
        x^* = x_0 + x_1\eta+x_2\eta^2 + \mathcal{O}(\eta^3),\label{eq:35}
    \end{align}
we have 
    \begin{align}
        y &\approx (1-(2+\delta_c)x_0) -(\delta_cx_0+(2+\delta_c)x_1)\eta \nonumber\\
        &\qquad - (\delta_cx_0+(2+\delta_c)x_2+\delta_cx_1)\eta^2 +\mathcal{O}(\eta^3).\label{eq:36}
    \end{align}
Substituting into (\ref{eq:32}) and equating coefficients to 0, we find an $\eta$-order expansion of the approximate equilibrium value $x^*$ as
    \begin{equation}
        x^* \approx \frac{1}{\lambda\sigma+\mu\delta_c +\mu -\delta_c}\eta+\mathcal{O}(\eta^2). \label{eq:37}
    \end{equation}
Using the relation $w^* = \frac{\sigma}{1-\eta}x^* = \sigma x^* +\mathcal{O}(\eta),$ we have
\begin{equation}
    w^* \approx \frac{\sigma}{\lambda\sigma+\mu\delta_c +\mu -\delta_c}\eta +\mathcal{O}(\eta^2). \label{eq:38}
\end{equation}
\begin{figure}
\centering
\caption{Bifurcation diagrams for the $\delta \approx \delta_c$ case}
\begin{subfigure}{.75\textwidth}
\centering
\includegraphics[scale=.4]{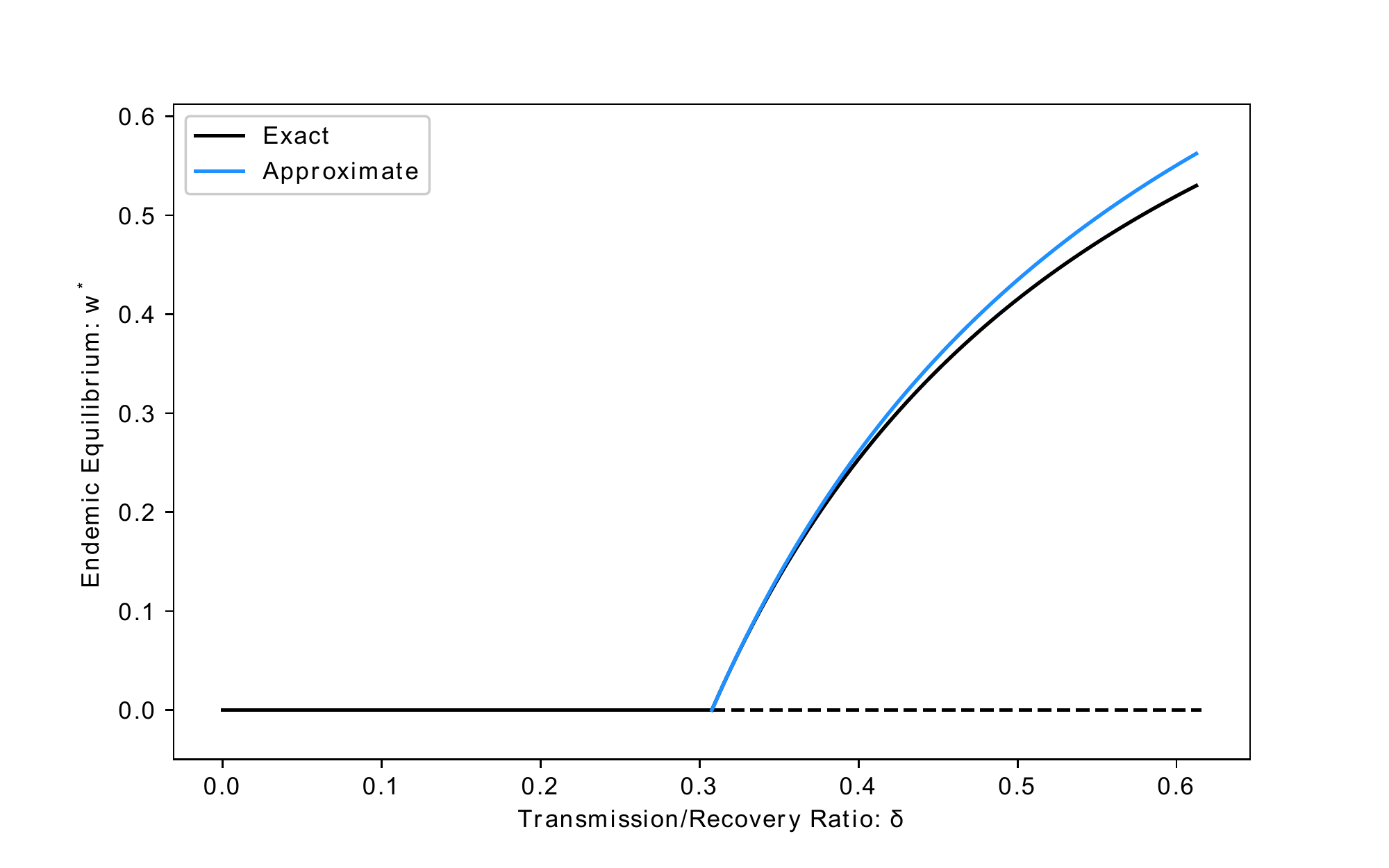}
\caption{Exact and approximate endemic equilibrium prevalence for a bimodal network with 5000 degree 3 nodes and 5000 degree 5 nodes. Moments of the degree distribution are $\kone = 4, \ktwo = 17, \kthree = 76$.}
\label{fig:1a}
\end{subfigure}
\begin{subfigure}{.75\textwidth}
\centering
\includegraphics[scale=.4]{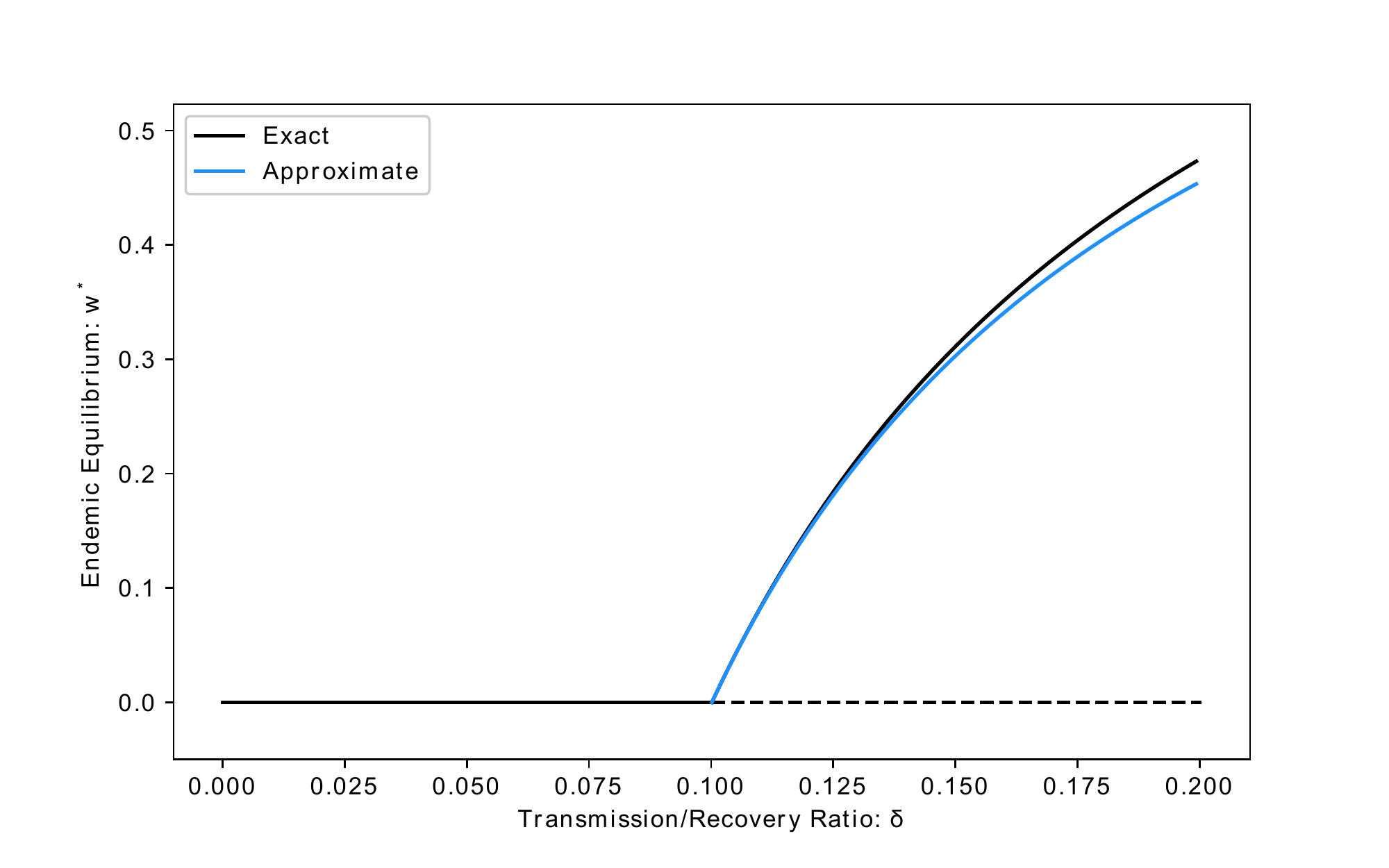}
\caption{Exact and approximate endemic equilibrium prevalence for a configuration-model network with a Poisson degree distribution with $\kone = 10$. Higher moments of the degree distribution are $\ktwo \approx 110, \kthree \approx 1309$.}
\label{fig:1b}
\end{subfigure}
\label{fig:1}
\end{figure}

To demonstrate the efficacy of this approximation, we compare the approximation (\ref{eq:38}) to the actual endemic equilibrium using bifurcation diagrams (Fig. \ref{fig:1}). We consider two example configuration model random networks \citep{molloy_critical_1995} with $N=10,000$. In Fig. \ref{fig:1a}, a bimodal network is considered with 5000 degree 3 nodes and 5000 degree 5 nodes. In Fig. \ref{fig:1b}, a network with a Poisson degree distribution (with average degree $\kone =10$) is considered. As is clear in both examples, the agreement between the actual and approximate endemic equilibrium is quite good near the epidemic threshold.

%%%%%%%%%%%%%%%%%%%%%%%%%%%%%%%%%%%%%%%%%%%%%%%%%%%%%%%%%%%%%%%%%%%%%%%%%%%%%%%%%%%   
\subsection{Case 2: Far away from the epidemic threshold ($\delta >> \delta_c$)}

For $\delta >> \delta_c,$  our small parameter of choice is $\epsilon = \delta_c/\delta$. We can rewrite (\ref{eq:29}) and (\ref{eq:30}) in terms of this parameter:
\begin{align}
    0 &= \epsilon^2(1-y-2x)(x+y)^2 \nonumber\\
    &\qquad - \epsilon\left(\delta_cx(x+y)^2+\lambda x^2 +\mu x^2(x+y)\right)+\lambda\sigma x^3, \label{eq:39}\\
    0 &= \epsilon^2(x+y)^2 -\epsilon\left(\lambda y +\mu y (x+y)\right) +\lambda\sigma xy. \label{eq:40}
\end{align}
When $\delta >> \delta_c,$ the transmission rate $\tau$ is large relative to the recovery rate $\gamma.$ Thus, we expect the disease to affect much of the population, and consequently there will be very few remaining $[SS]$ links, and therefore $y \approx 0.$

Solving $P(\phi,0)=0$ for $\phi$ yields
\begin{equation}
\phi(\epsilon) = \frac{\epsilon^2-\lambda\epsilon}{2\epsilon^2+(\delta_c+\mu)\epsilon-\lambda\sigma}\label{eq:41},    
\end{equation}    
and slope of the linearization is then
\begin{equation}
\psi(\epsilon)=-\frac{P_x(\phi,0)}{P_y(\phi,0)} = -\frac{(\epsilon-\lambda)(2\epsilon^2+(\delta_c+\mu)\epsilon-\lambda\sigma)}{\epsilon(\epsilon^2-(\mu+5\lambda)\epsilon-\lambda(2\delta_c+\mu-2\sigma))},  \label{eq:42} 
\end{equation}
so 
\begin{equation}
    y \approx \psi(x-\phi).\label{eq:43}
\end{equation}
Next, we seek to expand $y$ in terms of $\epsilon$ only. The relevant expansions for $\phi,\psi,$ and $x$ are:
\begin{align}
    \phi(\epsilon) &= \frac{1}{\sigma}\epsilon + \frac{\delta_c+\mu-\sigma}{\lambda\sigma^2}\epsilon^2 + \mathcal{O}(\epsilon^3), \label{eq:44}\\
    \psi(\epsilon) &= \frac{\lambda\sigma}{2\delta_c+\mu-2\sigma}\epsilon^{-1}-\frac{2\delta_c^2+3\delta_c\mu+\sigma(5\lambda+2\sigma)+\mu^2}{(2\delta_c+\mu-2\sigma)^2}+\mathcal{O}(\epsilon),\label{eq:45} \\
    x(\epsilon) &= x_0+x_1\epsilon+x_2\epsilon^2 + \mathcal{O}(\epsilon^2). \label{eq:46}
\end{align}

To ease the writing of coefficients, we let $\phi_\alpha$ and $\psi_\alpha$ refer to the coefficients on $\epsilon^\alpha$ for the respective series. From this, it follows that 
    \begin{align}
        y &\approx (\psi_{-1}x_0)\epsilon^{-1} + (\psi_{-1}x_1+\psi_0x_0-\psi_{-1}\phi_1) \nonumber \\ 
        &\qquad +(\psi_{-1}x_2+\psi_1x_0+\psi_0x_1-\psi_{-1}\phi_2-\psi_0\phi_1)\epsilon + \mathcal{O}(\epsilon^2). \label{eq:47}
    \end{align}
Substituting into (\ref{eq:40}), and equating the coefficients to 0, we find that we need the coefficients up to order $\epsilon^4$ in order to find a $\epsilon^2$ order expansion of the approximate equilibrium value of $x^*.$ The result is
    \begin{equation}
        x^* \approx \frac{1}{\sigma}\epsilon + \frac{\delta_c+\mu-\sigma}{\lambda\sigma^2}\epsilon^2+\mathcal{O}(\epsilon^3).\label{eq:48}
    \end{equation}
Finally, as $w^* =\sigma\epsilon^{-1}x^*,$ we arrive at an $\epsilon-$order approximation for size of the endemic steady state as 
    \begin{equation}
        w^* \approx 1 + \frac{\delta_c+\mu-\sigma}{\lambda\sigma}\epsilon +\mathcal{O}(\epsilon^2). \label{eq:49}
    \end{equation}

\begin{figure}
\centering
\caption{Bifurcation diagrams for the $\delta >> \delta_c$ case}
\begin{subfigure}{.75\textwidth}
\centering
\includegraphics[scale=.4]{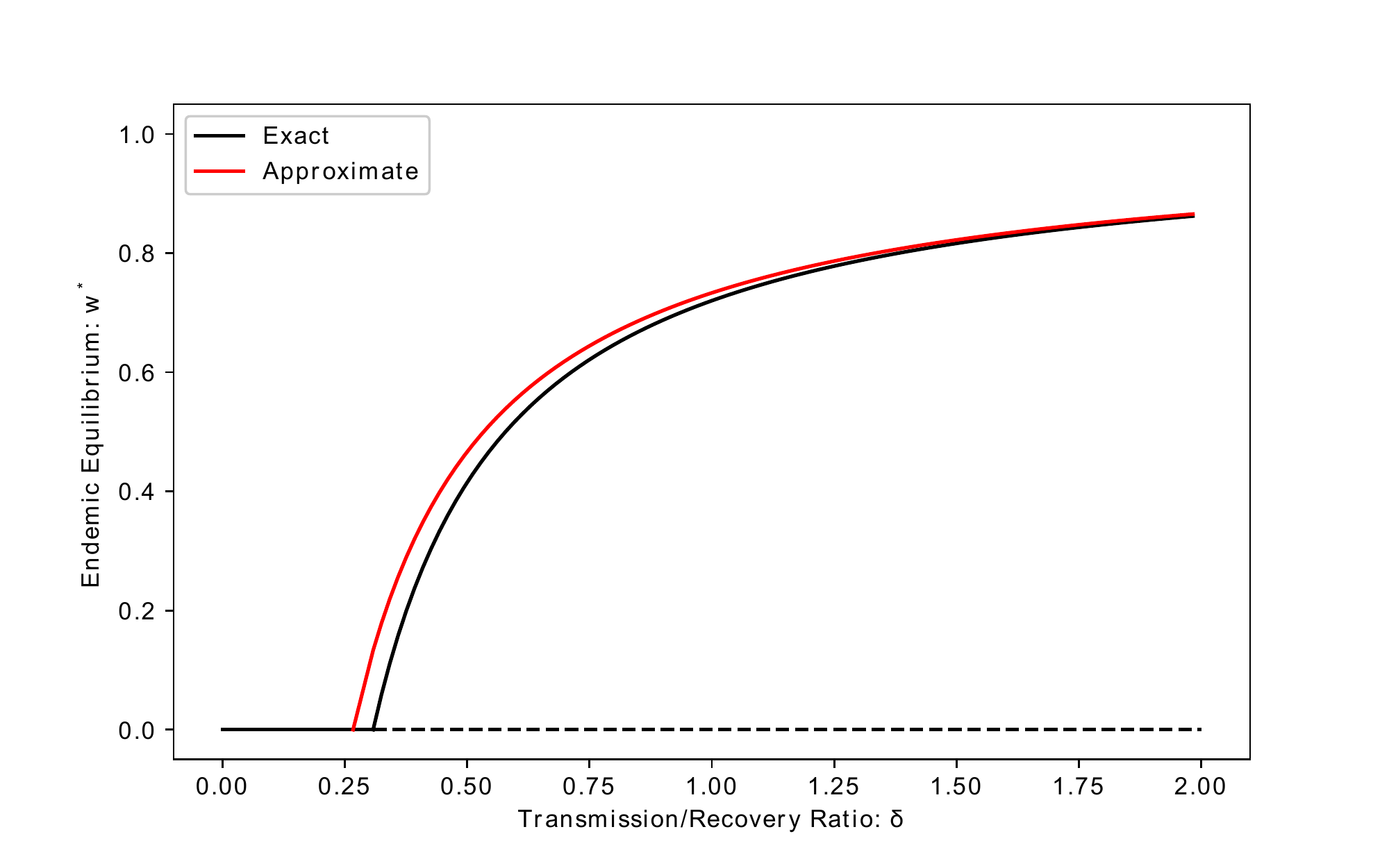}
\caption{Exact and approximate endemic equilibrium prevalence for a bimodal network with 5000 degree 3 nodes and 5000 degree 5 nodes. Moments of the degree distribution are $\kone = 4, \ktwo = 17, \kthree = 76$.}
\label{fig:2a}
\end{subfigure}
\begin{subfigure}{.75\textwidth}
\centering
\includegraphics[scale=.4]{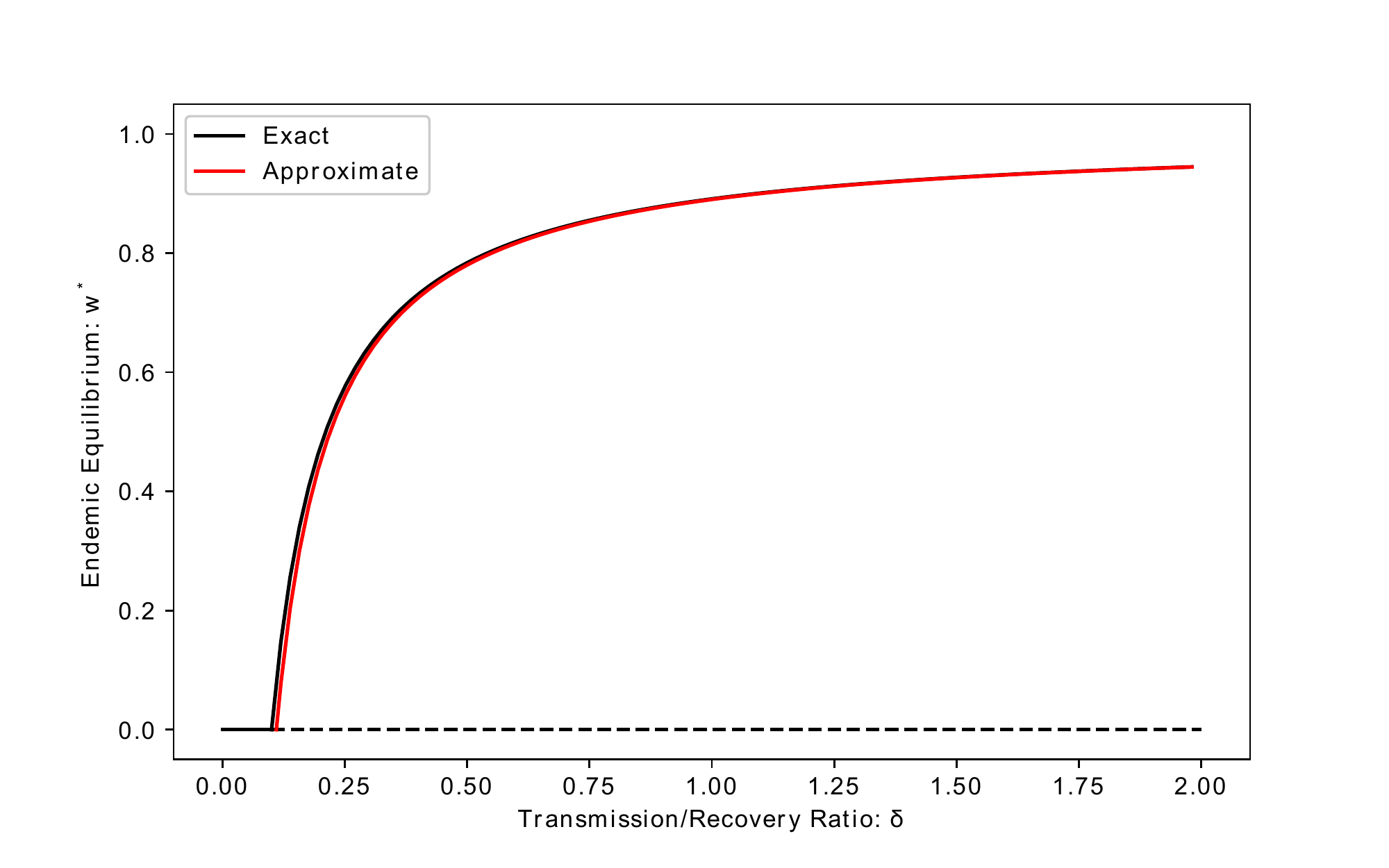}
\caption{Exact and approximate endemic equilibrium prevalence for a configuration-model network with a Poisson degree distribution with $\kone = 10$. Higher moments of the degree distribution are $\ktwo \approx 110, \kthree \approx 1309.$}
\label{fig:2b}
\end{subfigure}
\label{fig:2}
\end{figure}

\noindent As with the $\delta \approx \delta_c$ case, we compare the approximation (\ref{eq:49}) to the actual endemic equilibrium on the same networks as previously described. Again, the agreement is quite good, even for relatively small values of $\delta.$
% Moreover, this can be expressed entirely in terms of the original parameters and variables as
    
%     \begin{equation}
%         [I]^*/N = 1+\frac{\kthree-2\langle k ^2\rangle\kone+\kone^3}{\ktwo^2-\kthree\kone}\frac{\gamma}{\tau}.
%     \end{equation}

%%%%%%%%%%%%%%%%%%%%%%%%%%%%%%%%%%%%%%%%%%%%%%%%%%%%%%%%%%%%%%%%%%%%%%%%%%%%%%%%%%%
%%%%%%%%%%%%%%%%%%%%%%%%%%%%%%%%%%%%%%%%%%%%%%%%%%%%%%%%%%%%%%%%%%%%%%%%%%%%%%%%%%%
\subsection{Sensitivity Analysis}\label{sensitivity}

With any model of infectious disease, its implications in preventing or mitigating spread should be considered. For network models, some pharmaceutical and non-pharmaceutical interventions can alter the contact network structure in the effort to contain or mitigate outbreaks \citep{salathe_dynamics_2010}. For an SIS-type disease, particularly when containment is impossible, one such goal may be to decrease the size of the endemic equilibrium. To that end, we examine the sensitivity of our approximations of $w^*$ to network parameters in the SCPW model. One benefit of explicit asymptotic expressions for the endemic equilibrium is that sensitivity analyses are straightforward to implement.

For a fixed $\delta,$ we have a three-dimensional parameter space. To visualize these parameter combinations, we use two-dimensional heat maps taken at slices of the third network parameter. In this case, we have decided to look at several fixed values of $\kthree,$ and draw sensitivity heat maps in the variables $(\kone,\ktwo).$ Further complicating matters is the fact that moments of a distribution are subject to many inequalities which restrict the domain of the sensitivity heat maps. Two natural restrictions to include are the results of Jensen's Inequality and the Cauchy-Schwarz Inequality respectively:
\begin{align*}
    \ktwo &\geq \kone^2, \\
    \ktwo^2 &\leq \kthree \kone.
\end{align*}
For a fixed value of $\kthree,$ these restrictions give a wedge-shaped feasible region of $(\kone, \ktwo).$ We plot the sensitivities for $\kthree = 20, 100,$ and $400$ to display a range of possible parameter combinations. 

\begin{figure}
\caption{Sensitivities for the $\delta \approx \delta_c$ approximation}
\begin{subfigure}{\textwidth}
\centering
\includegraphics[scale=.25]{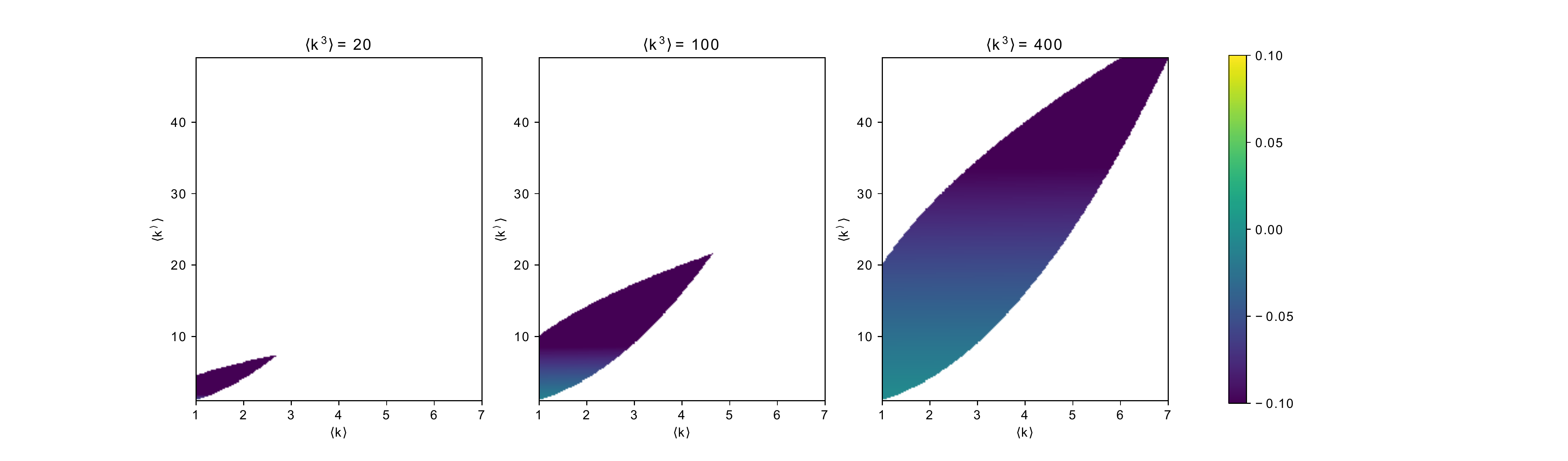}\label{fig:3a}
\caption{$\protect\frac{\partial w^*}{\partial \kone} $}
\end{subfigure}
\begin{subfigure}{\textwidth}
\centering
\includegraphics[scale=.25]{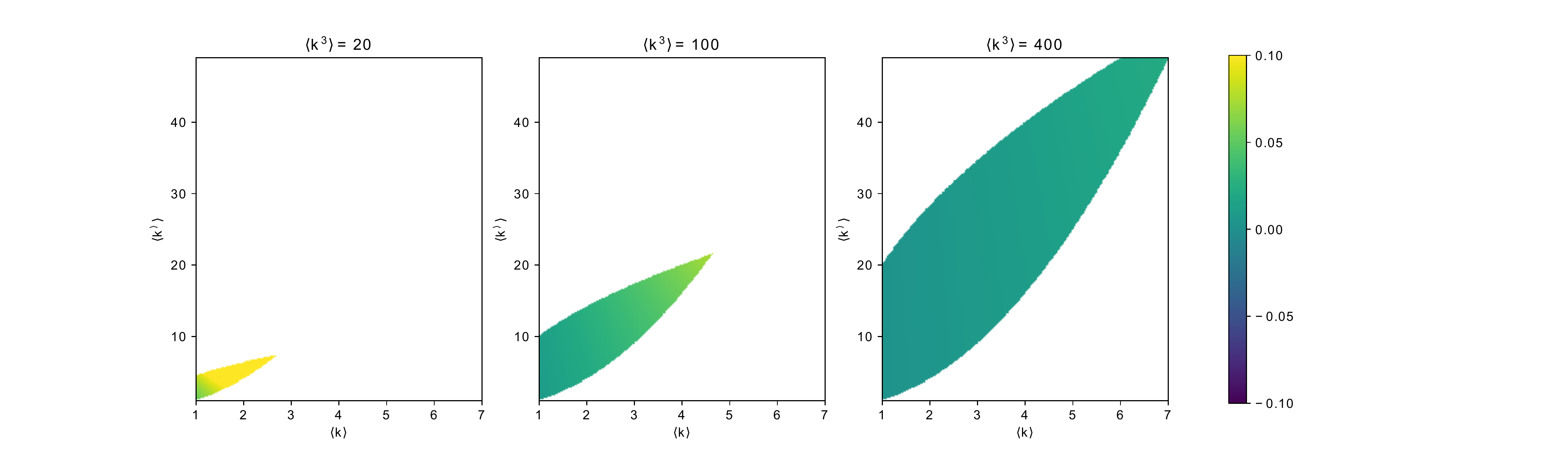}
\caption{$\protect\frac{\partial w^*}{\partial \ktwo}$} \label{fig:3b}
\end{subfigure}
% \begin{subfigure}{.8\textwidth}
% \includegraphics[scale=.4]{sens_nearz.pdf}
% \caption{$\protect\frac{\partial w^*}{\partial \kthree}$} \label{fig:1c}
%\end{subfigure}
\label{fig:3}
\end{figure}

In the $\delta\approx \delta_c$ case, calculating the partial derivatives is straightforward. To compute the sensitivities, we evaluate the partial derivatives at the epidemic threshold: $\delta = \delta_c.$ Table \ref{tab:1} shows the expressions for these sensitivities, and Fig. \ref{fig:3} shows corresponding plots. Clearly $\frac{\partial w^*}{\partial \kone} \leq 0$ and $\frac{\partial w^*}{\partial \kone} \geq 0,$ with more extreme values near the upper-right corner of the feasible region. 

\begin{table}
\begin{center}
\caption{Partial Derivatives for $\delta\approx \delta_c$}
\begin{tabular}{ r c l}
%\hline\\
%\multicolumn{3}{c}{Partial Derivatives for $\delta \approx \delta_c$}\\
\hline \\
$\left.\displaystyle\frac{\partial w^*}{\partial \kone }\right\vert_{\delta = \delta_c}$ & =  
&$\displaystyle-\frac{\ktwo}{\kone - 2\ktwo +\kthree}$\\
&&\\
$\left.\displaystyle\frac{\partial w^*}{\partial \ktwo }\right\vert_{\delta = \delta_c}$ & = & $\displaystyle\frac{\kone}{\kone - 2\ktwo +\kthree}$\\
&&\\
$\left.\displaystyle\frac{\partial w^*}{\partial \kthree }\right\vert_{\delta=\delta_c}$ & = & 0\\
&&\\
\hline
\end{tabular}
\label{tab:1}
\end{center}
\end{table}

\begin{figure}
\caption{Sensitivities for the $\delta >> \delta_c$ approximation}
\begin{subfigure}{\textwidth}
\centering
\includegraphics[scale=.25]{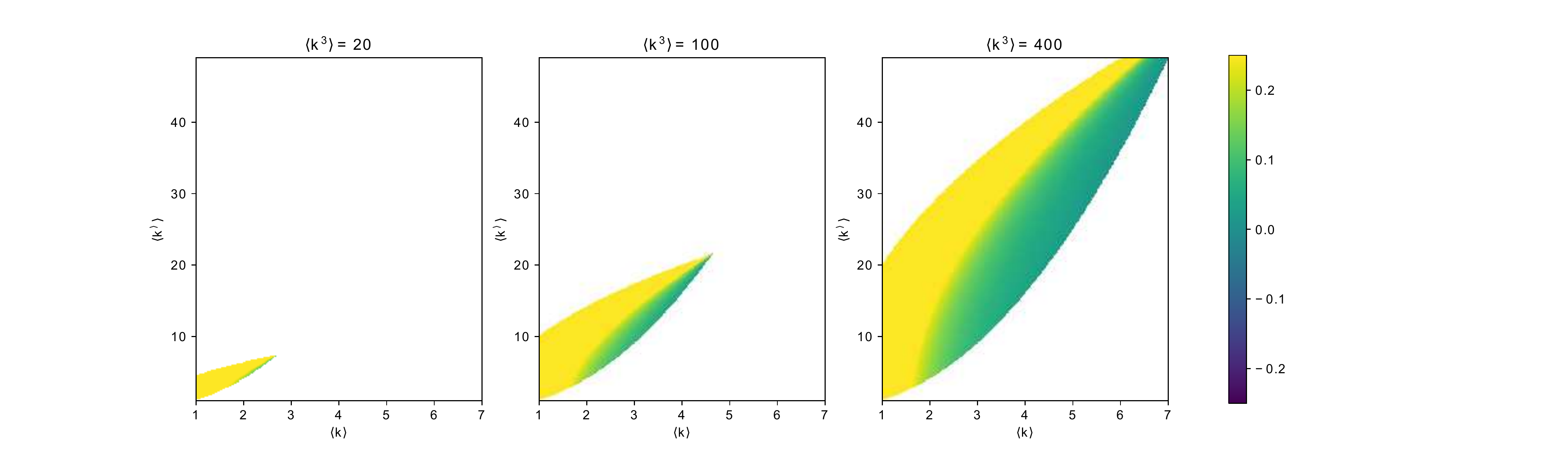}\label{fig:4a}
\caption{$\protect\frac{\partial w^*}{\partial \kone} $}
\end{subfigure}
\begin{subfigure}{\textwidth}
\centering
\includegraphics[scale=.25]{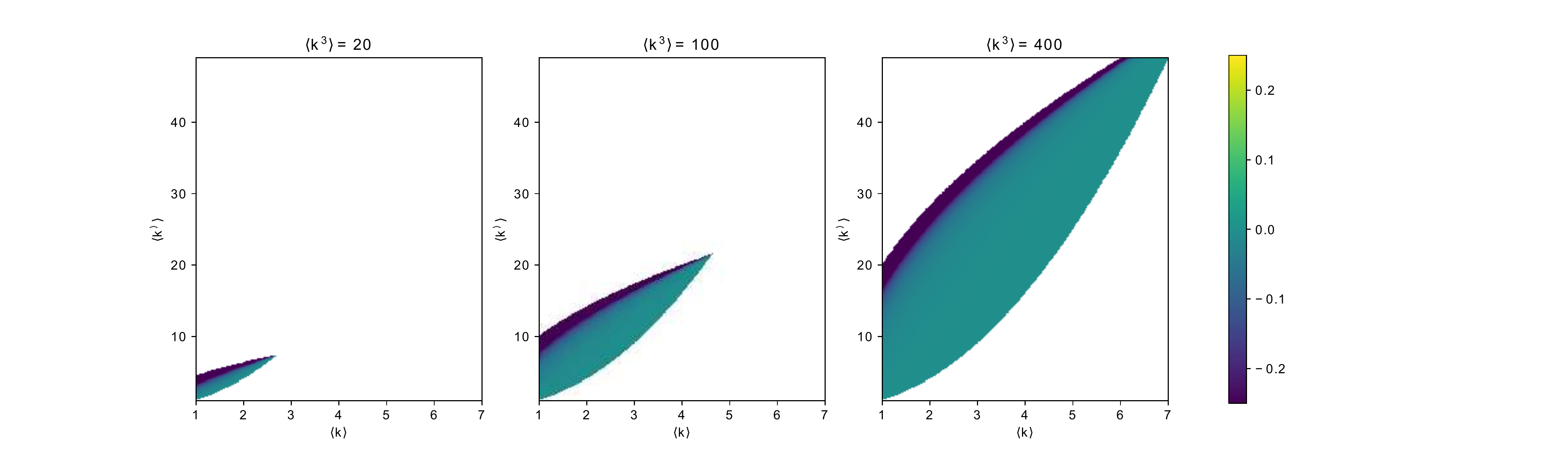}
\caption{$\protect\frac{\partial w^*}{\partial \ktwo}$} \label{fig:4b}
\end{subfigure}
\begin{subfigure}{\textwidth}
\centering
\includegraphics[scale=.25]{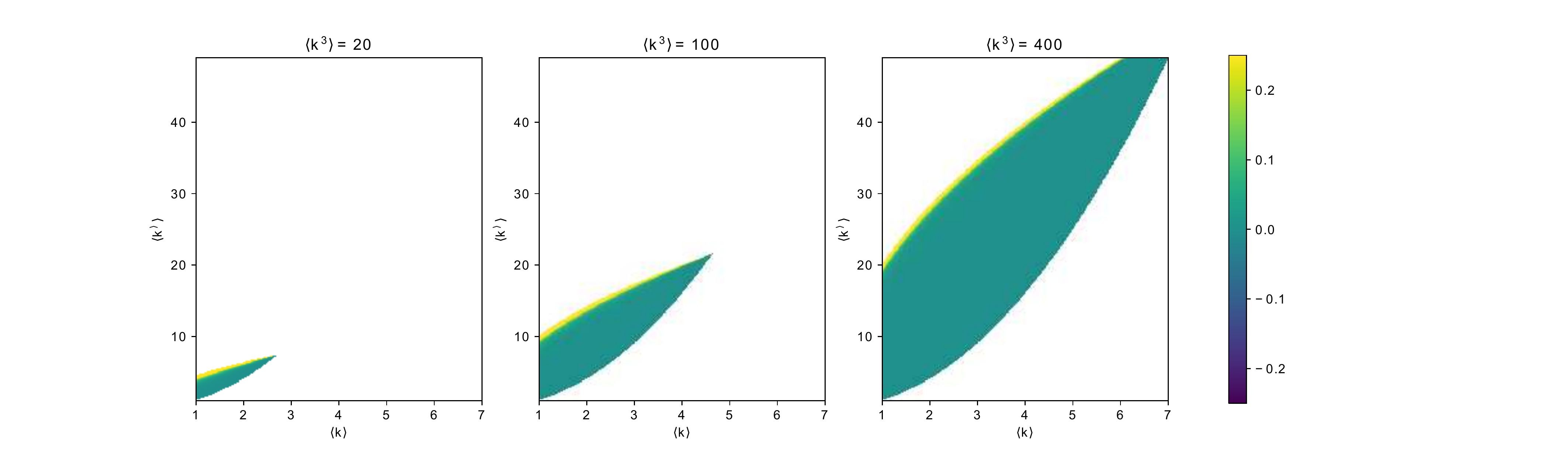}
\caption{$\protect\frac{\partial w^*}{\partial \kthree}$} \label{fig:4c}
\end{subfigure}
\label{fig:4}
\end{figure}

For the $\delta >> \delta_c$ case, the partial derivatives (Table \ref{tab:2}) all depend on a factor of $1/\delta,$ so the choice of $\delta$ for computing sensitivities does not affect the relative magnitudes of the partial derivatives. For convenience, we select $\delta=1.5.$ The sensitivity plots in Fig. \ref{fig:4} show that $\frac{\partial w^*}{\partial \kone} \geq 0, \frac{\partial w^*}{\partial \ktwo} \leq 0,$ and $\frac{\partial w^*}{\partial \kthree} \geq 0,$ with the greatest sensitivity near the curve $\ktwo^2 = \kthree \kone,$ though the large magnitude appears to be due to the partial derivatives being undefined there.

\begin{table}
\begin{center}
\caption{Partial Derivatives for $\delta >>\delta_c$}
\begin{tabular}{ r c l}
%\hline\\
%\multicolumn{3}{c}{Partial Derivatives for $\delta \approx \delta_c$}\\
\hline \\
$\displaystyle\frac{\partial w^*}{\partial \kone }$ & =  
&$\displaystyle\frac{\kthree^2+3\kone^2\ktwo^2-2(\kone^3\kthree+\ktwo^3)}{(\ktwo^2 - \kthree\kone)^2}\frac{1}{\delta}$\\
&&\\
$\displaystyle\frac{\partial w^*}{\partial \ktwo }$ & = & $\displaystyle-\frac{2(\kone^2-\ktwo)(\kone\ktwo - \kthree)}{(\ktwo^2 - \kthree\kone)^2}\frac{1}{\delta}$\\
&&\\
$\displaystyle\frac{\partial w^*}{\partial \kthree }$ & = & $\displaystyle\frac{(\kone^2-\ktwo)^2}{(\ktwo^2 - \kthree\kone)^2}\frac{1}{\delta}$\\
&&\\
\hline
\end{tabular}
\label{tab:2}
\end{center}
\end{table}

A significant observation from these sensitivities is that $\frac{\partial w^*}{\partial \kone}$ and $\frac{\partial w^*}{\partial \ktwo}$ change signs depending on the regime considered. If the goal on an intervention is to reduce the size of the endemic equilibrium, near the epidemic threshold, this can be accomplished in principle by increasing $\kone$ or decreasing $\ktwo,$ which will in effect increase $\delta_c$ as well. This is intuitive, as an effort to push the system below the epidemic threshold would also decrease the endemic equilibrium for a fixed $\delta$. However, in the $\delta >> \delta_c$ regime, the system is far from the epidemic threshold, and reducing the size of the endemic equilibrium can be accomplished by decreasing $\kone$ or increasing $\ktwo.$ This suggests that containment and mitigation strategies that depend on altering the structure of the contact network may require different goals in terms of the moments of the degree distribution.
%%%%%%%%%%%%%%%%%%%%%%%%%%%%%%%%%%%%%%%%%%%%%%%%%%%%%%%%%%%%%%%%%%%%%%%%%%%%%%%%%%
%%%%%%%%%%%%%%%%%%%%%%%%%%%%%%%%%%%%%%%%%%%%%%%%%%%%%%%%%%%%%%%%%%%%%%%%%%%%%%%%%%

\section{Conclusion}\label{Conclusion}

In this paper, we have analyzed the super compact pairwise model presented in \citet{simon_super_2016-1}. A non-dimensional version of the model was considered, and a bifurcation analysis was performed demonstrating that the SCPW and CPW models share an epidemic threshold. Moreover, we derived approximate formulas for the endemic equilibrium in two regimes: when the transmission/recovery ratio is near the epidemic threshold, and far away from it. While the asymptotic techniques used here are ad hoc, similar techniques may prove fruitful in other low-dimensional models of infectious disease spread on networks. However, an exact expression for the endemic equilibrium remains elusive. 

Before explaining the advantages of our approach, we acknowledge two limitations of our approximation. First, approximations of the endemic equilibrium for diseases between the two regimes is lacking. Second, while the examples of simulated networks show good agreement between the exact and approximate prevalence, we have not quantified the approximation error generally. As such, there may be types of networks for which our approximation of the endemic equilibrium is less accurate or inappropriate. 

Our approximation of the endemic equilibrium is very useful in providing a more detailed look into the interactions of the moments of the degree distribution as they relate to the size of an outbreak. This has implication for disease control measures, particularly those that work by altering the contact network structure. Our results suggest that for SIS-type diseases, strategies to contain (near the epidemic threshold) or mitigate (far away from the epidemic threshold) an outbreak may require different goals. In the mitigation scenario where the prevalence is high, measures might be employed that decrease the first moment $\kone$ of the degree distribution. In effect, this may mean initiatives aimed at reducing the number of contacts of individuals alone. On the other hand, in the containment scenario where the prevalence is low, decreasing the second moment $\ktwo$ may be efficient. When couched in degree distribution terms this goal is hard to conceptualize, but using probability generation functions \citep{newman_random_2001} one can show that $\ktwo$ is the average number of first and second neighbors of nodes in the network. Thus, measures that reduce both the contacts of individuals and their partners are effective in this scenario. This suggests the importance of contact tracing.

Our results complement the findings of \citet{eames_modeling_2002}, who observed that the effectiveness of two common control measures, screening and contact tracing, depend on the prevalence at the endemic equilibrium. Screening, which targets and treats individuals, is efficient when the prevalence is high. Contact tracing, which targets and treats individuals and their partners, if efficient when the prevalence is low. Unlike this paper, Eames and Keeling implement these measures through epidemiological parameters (rather than through changing network structure). In this way, our results can be viewed as a network-structure analog for their conclusions and confirm that control measures appropriate in a network setting can be found. Further work in this area may include investigating this phenomenon with alternative models of SIS diseases on networks.

%%%%%%%%%%%%%%%%%%%%%%%%%%%%%%%%%%%%%%%%%%%%%%%%%%%%%%%%%%%%%%%%%%%%%%%%%%%%%%%%%%
%%%%%%%%%%%%%%%%%%%%%%%%%%%%%%%%%%%%%%%%%%%%%%%%%%%%%%%%%%%%%%%%%%%%%%%%%%%%%%%%%%
\section*{Acknowledgment}  This research supported in part by NSF Grant DMS - 1817124 to AH.

\begin{appendices}
\numberwithin{equation}{section}

\section{Bifurcation and Endemic Steady State} \label{AppendixA}

We begin with Theorem 4.1 from \citet{castillo-chavez_dynamical_2004}, referring to the specific conditions that will be relevant for this analysis. Consider a system of ODEs with a parameter $\phi:$
\begin{equation} 
\frac{dx}{dt} = \mathbf{F}(x,\phi),\:\:\: \mathbf{F}:\mathbb{R}^n\times\mathbb{R}\to \mathbb{R}^n \text{ and } \mathbf{F}\in C^2(\mathbb{R}^n\times\mathbb{R}). \label{eq:A1} 
\end{equation}
Assume that $0$ is an equilibrium for all values of $\phi.$ Assume further that $D_xf(0,0) = \left(\frac{\partial F_i}{\partial x_j}(0,0)\right)$ is the linearization matrix of (\ref{eq:A1}) around the equilibrium 0 and with $\phi=0$, and zero is a simple eigenvalue of this matrix with all other eigenvalues having negative real parts. Assume as well that this matrix has a nonnnegative right eigenvector $\mathbf{w}$ and left eigenvector $\mathbf{v}$ corresponding to the zero eigenvalue. Let $F_k$ be the $k$th component of $f$ and 
\begin{align}
    a &= \sum_{k,i,j=1}^n v_kw_iw_j\frac{\partial^2
    F_k}{\partial x_i\partial x_j}(0,0), \label{eq:A2}\\
    b &= \sum_{k,i=1}^n v_k w_i \frac{\partial^2 F_k}{\partial x_i\partial \phi}(0,0).\label{eq:A3}
\end{align}
If $a<0$ and $b>0,$ then when $\phi$ changes from negative to positive, 0 changes its stability from stable to unstable. Correspondingly, a negative unstable equilibrium becomes positive and locally asymptotically stable.

We apply this theorem to $(\ref{eq:17})-(\ref{eq:19}),$ where the equilibrium occurs at $w=x=z=0.$ Moreover, we define a bifurcation parameter $\phi = \delta-\delta_c,$ so $\phi = 0$ corresponds to $\delta=\delta_c,$ and $\frac{\partial}{\partial \phi} = \frac{\partial}{\partial \delta}.$ For consistency with previously established notation, we will treat $\delta$ as our parameter, with $\phi$ increasing through $0$ as $\delta$ increases through $\delta_c$. The Jacobian given in (\ref{eq:21}) when $w=0,x=0,z=0,$ and $\delta = \delta_c$ is 
\begin{equation}
J = \begin{bmatrix}-1 & \kone \delta_c & 0\\
0 & -\delta_c & 1\\ 0 &2\delta_c&-2\\ \end{bmatrix}.\label{eq:A4}
\end{equation}
and the characteristic polynomial is given by
\begin{equation}
0 = \xi(\xi+1)(\xi - (-2-\delta_c)).\label{eq:A5}
\end{equation}
The left and right eigenvectors ($\mathbf{v}$ and $\mathbf{w}$ respectively) corresponding to the eigenvalue $\xi =0$ are 
\begin{equation}
\mathbf{v} = \begin{bmatrix} 0 & 2 & 1\end{bmatrix}, \mathbf{w} = \begin{bmatrix} \kone & \delta_c^{-1}&1\end{bmatrix}^T. \label{eq:A6}
\end{equation}
To compute $a$ and $b$, it is convenient to express (\ref{eq:A2}) and (\ref{eq:A3}) in matrix-vector form:
\begin{align}
    a &= \mathbf{w}^T (2H_2 + H_3)\mathbf{w},\label{eq:A7}\\
    b &= \mathbf{v}\frac{\partial J}{\partial \delta}(0,\delta_c) \mathbf{w},\label{eq:A8}
\end{align}
where $H_2$ and $H_3$ are the Hessians of $F_2$ and $F_3$ respectively at $\Vec{0}$. These Hessians are 
\begin{equation}
H_2 = \begin{bmatrix} 0 & -\dfrac{\alpha\delta_c}{\kone}&0\\ -\dfrac{\alpha\delta_c}{\kone} & -2-2\beta\delta_c & \dfrac{\alpha\delta_c}{\kone}\\ 0 & \dfrac{\alpha\delta_c}{\kone}&0\\\end{bmatrix}, H_3 = \begin{bmatrix} 0&0&0\\0&4&0\\0&0&0 \end{bmatrix}. \label{eq:A9}
\end{equation}
Thus,
\begin{align}
    a &= \begin{bmatrix}\kone & \delta_c^{-1} & 1 \end{bmatrix}\begin{bmatrix} 0 & -\dfrac{2\alpha\delta_c}{\kone}&0\\ -\dfrac{2\alpha\delta_c}{\kone} & -4\beta\delta_c & \dfrac{2\alpha\delta_c}{\kone}\\ 0 & \dfrac{2\alpha\delta_c}{\kone}&0\\\end{bmatrix}\begin{bmatrix}\kone \\ \delta_c^{-1} \\ 1 \end{bmatrix}\nonumber \\
    &=\begin{bmatrix}\kone & \delta_c^{-1} & 1 \end{bmatrix}\begin{bmatrix} -\dfrac{2\alpha}{\kone}\\ -2\alpha\delta_c-4\beta + \dfrac{2\alpha\delta_c}{\kone}\\ \dfrac{2\alpha}{\kone}\\\end{bmatrix}\nonumber \\
    &= -2\alpha -2\alpha-4\beta/\delta_c+2\dfrac{\alpha}{\kone}+2\frac{\alpha}{\kone}\nonumber \\
    &= -4\left(\alpha\left(\frac{1}{\kone} +1\right)+\beta\left(\frac{\ktwo}{\kone}+1\right) \right)\nonumber \\
    &=-4\left(\frac{\kthree}{\kone}-1\right). \label{eq:A10}
\end{align}
As $\kthree> \kone,$ it follows that $a<0.$

The computation for $b$ is simpler. We note that 
\begin{equation}
    \frac{\partial J}{\partial \delta}(0,\delta_c) = \begin{bmatrix} 0 & \kone & 0 \\ 0 & \delta_c^{-1}-1 & 0 \\ 0 & 2& 0\\ \end{bmatrix}.\label{eq:A11}
\end{equation}
Thus
\begin{align}
    b &= \begin{bmatrix} 0 & 2 & 1 \end{bmatrix}\begin{bmatrix} 0 & \kone & 0 \\ 0 & \delta_c^{-1}-1 & 0 \\ 0 & 2& 0\\ \end{bmatrix}\begin{bmatrix}\kone \\ \delta_c^{-1} \\ 1 \end{bmatrix}\nonumber \\
    &=\begin{bmatrix} 0 & 2 & 1 \end{bmatrix}\begin{bmatrix} 0 & \kone \delta_c^{-1}& 0 \\ 0 & \delta_c^{-1}(\delta_c^{-1}-1) & 0 \\ 0 & 2\delta_c^{-1}& 0\\ \end{bmatrix} \nonumber\\
    &= 2\delta_c^{-1}(\delta_c^{-1}-1)+2\delta_c^{-1} = 2\delta_c^{-2} >0. \label{eq:A12}
\end{align}
Finally, as $a<0$ and $b>0,$ we conclude that as $\delta$ increases through $\delta_c,$ a positive, asymptotically stable equilibrium emerges, which is the endemic equilibrium.

%%%%%%%%%%%%%%%%%%%%%%%%%%%%%%%%%%%%%%%%%%%%%%%%%%%%%%%%%%%%%%%%%%%%%%%%%%%%%%%%%%
%%%%%%%%%%%%%%%%%%%%%%%%%%%%%%%%%%%%%%%%%%%%%%%%%%%%%%%%%%%%%%%%%%%%%%%%%%%%%%%%%%

\section{Asymptotic Approximations of the Endemic Equilibrium}

The full derivations of the approximations (\ref{eq:38}) and (\ref{eq:49}) are presented in this appendix. 
%%%%%%%%%%%%%%%%%%%%%%%%%%%%% %%%%%%%%%%%%%%%%%%%%%%%%%%%%%%%%%%%%%%%%%%%%%%%%%%%%%%
\subsection{Near the epidemic threshold ($\delta \approx \delta_c$)}

We begin with (\ref{eq:31}) and (\ref{eq:32}) and seek the linear approximation of $P(x,y) = 0$ at $(0,1).$ We compute
\begin{align}
    \frac{\partial P}{\partial x} &= 2(1-\eta)^2\left((1-y-2x)(x+y)-(x+y)^2\right) \nonumber\\ 
    &\qquad-(1-\eta)\left(\delta_c (x+y)(3x+y)+2\lambda x +\mu x (3x+2y)\right) + 3\lambda \sigma x^2,\label{eq:B1}\\
    \frac{\partial P}{\partial y} &= (1-\eta)\left(-2\delta_c x(x+y)-\mu x^2+(1-\eta)(x+y)(2-5x-3y)\right).\label{eq:B2}
\end{align}
The slope of the linear approximation is then
\begin{equation}
    \left.-\frac{\partial P/\partial x}{\partial P/\partial y}\right\vert_{(0,1)} = -\frac{-2(1-\eta)^2-\delta_c(1-\eta)}{-(1-\eta)^2} = -2-\frac{\delta_c}{1-\eta},\label{eq:B3}
\end{equation}
and thus we approximate
\begin{equation}
    y \approx 1+\left(-2-\frac{\delta_c}{1-\eta}\right)x. \label{eq:B4}
\end{equation}

We now expand $x$ as $x = x_0 + x_1\eta+\dots$ and $\frac{\delta_c}{1-\eta} = \delta_c(1+\eta+\eta^2+\dots)$ as a geometric series. Incorporating these with (\ref{eq:B4}), we get the approximate expansion of $y$ as 
\begin{align}
    y &\approx 1 - \left(2+\delta_c(1+\eta+\eta^2+\dots )\right)(x_0 + x_1\eta + x_2\eta^2\dots)\nonumber\\
    &= 1-(2+\delta_c)x_0 - (\delta_c x_0+(2+\delta_c)x_1)\eta\nonumber\\ &\qquad-(\delta_c x_0 + \delta_c x_1 + (2+\delta_c)x_2)\eta^2+\dots \label{eq:B5}
\end{align}
For easier bookkeeping, define $y_\alpha$ to be the coefficient of $\eta^\alpha$ in (\ref{eq:B5}). As well, the following expansions will prove useful:
\begin{align}
    x^2 &= x_0^2 + 2x_0x_1\eta +(x_1^2+2x_0x_2)\eta^2+\dots, \label{eq:B6}\\
    y^2 &= y_0^2 + 2y_0y_1\eta +(y_1^2+2y_0y_2)\eta^2+\dots, \label{eq:B7}\\
    xy &= x_0y_0 + (x_0y_1+x_1y_0)\eta +(x_0y_2+x_1y_1+x_2y_0)\eta^2+\dots\label{eq:B8}
\end{align}
Now, we apply (\ref{eq:B5})-(\ref{eq:B8}) to (\ref{eq:32}) yielding
\begin{align}
0 &= (1-2\eta+\eta^2)\left(x_0^2 +2x_0y_0+y_0^2+2(x_0x_1+x_0y_1+x_1y_0 + y_0y_1)\eta + \dots \right)\nonumber \\
&-(1-\eta)\left(\lambda y_0+\mu y_0(x_0+y_0)+(\lambda y_1 +\mu(x_0y_1+x_1y_0+2y_0y_1))\eta+\dots\right)\nonumber\\
&+\lambda\sigma\left(x_0y_0 + (x_0y_1+x_1y_0)\eta +\dots \right).\label{eq:B9}
\end{align}
Equating the $\mathcal{O}(1)$ terms to zero, we have
\begin{align}
    0 &= x_0^2+2x_0y_0 +y_0^2-\lambda y_0 -\mu y_0(x_0+y_0)+\lambda\sigma x_0y_0 \nonumber\\
    &= (1-(1+\delta_c)x_0)^2-\lambda(1-(2+\delta_c)x_0)\nonumber\\
    &\qquad-\mu(1-(2+\delta_c)x_0)(1-(1-\delta_c)x_0) \nonumber\\
    &\qquad\qquad +\lambda\sigma x_0(1-(2+\delta_c)x_0) \nonumber\\
    &=1-2(1+\delta_c)x_0+x_0^2-\lambda +\lambda(2+\delta_c)x_0-\mu(1-(3+2\delta_c)x_0) \nonumber\\
    &\qquad\qquad -\mu(1+\delta_c)(2+\delta_c)x_0^2+\lambda\sigma x_0-\lambda\sigma(2+\delta_c)x_0^2 \nonumber\\
    &=\left(1-\lambda-\mu\right)+\left(\lambda\sigma + \lambda(2+\delta_c)+\mu(3+2\delta_c)-2(1+\delta_c)\right)x_0 \nonumber\\
    &\qquad\qquad+\left(1-\mu(1+\delta_c)(2+\delta_c)-\lambda\sigma(2+\delta_c)\right)x_0^2 \nonumber\\
    &= x_0\left[\lambda\sigma + \lambda(2+\delta_c)+\mu(3+2\delta_c)-2(1+\delta_c)\right. \nonumber\\
    &\qquad+ \left.\left(1-\mu(1+\delta_c)(2+\delta_c)-\lambda\sigma(2+\delta_c)\right)x_0\right].\label{eq:B10}
\end{align}
where we avail ourselves of (\ref{eq:22}) for the last equality. For the solution were interested, we have $x_0 = 0$ and $y_0 = 1.$

We rewrite (\ref{eq:B9}) as
\begin{align}
0 &= (1-2\eta+\eta^2)\left(1+(x_1 + 2y_1)\eta + \dots \right)\nonumber \\
&-(1-\eta)\left(\lambda+\mu+(\lambda y_1 +\mu(x_1+2y_1))\eta+\dots\right)\nonumber\\
&+\lambda\sigma\left(x_1\eta +\dots \right)\label{eq:B11}.
\end{align}
Equating the coefficients of the $\mathcal{O}(\eta)$ terms to zero gives
\begin{align}
    0 &= -2+2x_1+2y_1+(\lambda+\mu)-(\lambda y_1+\mu(x_1+2y_1))+\lambda\sigma x_1 \nonumber\\
    &= -2+2x_1-2(2+\delta_c)x_1+1+\lambda(2+\delta_c)x_1 \nonumber\\
    &\qquad-\mu(x_1-2(2+\delta_c)x_1)+\lambda\sigma x_1 \nonumber\\
    &= -1 +x_1\left(2-2(2+\delta_c)+\lambda(2+\delta_c)-\mu(1-2(2+\delta_c))+\lambda\sigma\right) \nonumber\\
    &= -1 + x_1\left(\lambda\sigma+\mu\delta_c+\mu-\delta_c\right). \label{eq:B12}
\end{align}
Thus,
\begin{equation}
    x_1 = \frac{1}{\lambda\sigma+\mu\delta_c+\mu-\delta_c}.\label{eq:B13}
\end{equation}
Now that we have a first order approximation of $x,$ we obtain an first order approximation of the endemic equilibrium:
\begin{align}
    w^* &= \frac{\sigma}{1-\eta}x^* \nonumber \\
    &= \sigma(1+\eta+\eta^2+\dots)(x_0+x_1\eta+\dots) \nonumber\\
    &= \sigma x_1 \eta + \mathcal{O}(\eta^2). \label{eq:B14}
\end{align}
and thus
\begin{equation}
    w^* \approx \frac{\sigma}{\lambda\sigma+\mu\delta_c+\mu-\delta_c}\eta + \mathcal{O}(\eta^2). \label{eq:B15}
\end{equation}

%%%%%%%%%%%%%%%%%%%%%%%%%%%%%%%%%%%%%%%%%%%%%%%%%%%%%%%%%%%%%%%%%%%%%%%%%%%%%%%%%%
\subsection{Far away from the epidemic threshold ($\delta >> \delta_c$)}

We begin with (\ref{eq:39}) and (\ref{eq:40}) and seek the linear approximation of $P(x,y) = 0$ at $(\phi,0)$ where $\phi$ is given by (\ref{eq:41}). We compute
\begin{align}
    \frac{\partial P}{\partial x} &= -2\epsilon^2(x+y)(3x+2y-1) \nonumber\\
    &\:\:-\epsilon\left(\delta_c(3x^2+4xy+y^2)+2\lambda x +\mu x(3x+2y)\right)+2\lambda\sigma x^2, \label{eq:B16}\\
    \frac{\partial P}{\partial y} &= \epsilon\left(-2\delta_c x(x+y)-\mu x^2-\epsilon(x+y)(5x+3y-2)\right). \label{eq:B17}
\end{align}
The slope of the linear approximation is then
\begin{align}
    \left.-\frac{\partial P/\partial x}{\partial P/\partial y}\right\vert_{(\phi,0)} &= \frac{-2\epsilon^2\phi(3\phi-1)-3\epsilon \delta_c\phi^2+3\lambda \phi + 3\mu \phi^2+2\lambda\sigma\phi^2}{\epsilon\left(-2\delta_c\phi^2-\mu\phi^2-\epsilon\phi(5\phi-2)\right)}\nonumber\\
    &= \frac{-(\epsilon^2-\epsilon\lambda)\phi}{\phi(-\epsilon \phi(2\delta_c+\mu)-\epsilon^2(5\phi-2))}\nonumber\\
    &= \frac{-(\epsilon-\lambda)(2\epsilon^2+\epsilon(\delta_c+\mu)-\lambda\sigma)}{\epsilon\left(\epsilon^2-\epsilon(\mu+5\lambda)+2\lambda \sigma-\lambda(2\delta_c+\mu)\right)}=\psi(\epsilon).\label{eq:B18}
\end{align}
Thus, the linear approximation at $(\phi,0)$ is 
\begin{equation}
    y \approx \psi(x-\phi). \label{eq:B19}
\end{equation}
We now expand $\psi$ and $\phi$ in powers of $\epsilon:$
\begin{align}
    \psi(\epsilon) &= \psi_{-1}\epsilon^{-1}+\psi_0 + \mathcal{O}(\epsilon)\nonumber\\
    &=\frac{\lambda\sigma}{2\delta_c+\mu-2\sigma}\epsilon^{-1}-\frac{2\delta_c^2+3\delta_c\mu+\sigma(5\lambda+2\sigma)+\mu^2}{(2\delta_c+\mu-2\sigma)^2}+\mathcal{O}(\epsilon), \label{eq:B20}\\
    \phi(\epsilon) &= \phi_1\epsilon + \phi_2\epsilon + \mathcal{O}(\epsilon^3)\nonumber\\
    &= \frac{1}{\sigma}\epsilon +\frac{\delta_c+\mu-\sigma}{\lambda\sigma^2}\epsilon^2+\mathcal{O}(\epsilon^3). \label{eq:B21}
\end{align}
Now, we expand $x$ as well and reorganize to express $y$ as a power series in $\epsilon:$
\begin{align}
    y &\approx \left(\psi_{-1}\epsilon^{-1}+\psi_0 +\dots \right)\left(\left(x_0+x_1\epsilon_1+x_2\epsilon^2+\dots\right)-\left(\phi_1\epsilon+\phi_2\epsilon^2+\dots\right)\right)\nonumber\\
    &= \psi_{-1}x_0\epsilon^{-1}+(\psi_{-1}x_1+\psi_0x_0-\psi_{-1}\phi_1)\nonumber\\
    &\qquad +\left(\psi_{-1}x_2+\psi_0x_1+\psi_1 x_0-(\psi_{-1}\phi_2+\psi_0 \phi_1)\right)\epsilon + \mathcal{O}(\epsilon^2).\label{eq:B22}
\end{align}
For easier bookkeeping, we define $y_\alpha$ to be the coefficient of $\epsilon^\alpha$ in (\ref{eq:B22}). Again, the following expansions will prove useful:
\begin{align}
    x^2 &= x_0^2 + 2x_0x_1\epsilon +(x_1^2+2x_0x_2)\epsilon^2+\dots,\label{eq:B23}\\
    y^2 &= y_{-1}\epsilon^{-2} + 2y_{-1}y_0\epsilon^{-1} +(y_0^2+2y_{-1}y_1)+\dots,\label{eq:B24}\\
    xy &= x_0y_{-1}\epsilon^{-1} + (x_0y_0+x_1y_{-1}) +(x_0y_1+x_1y_0+x_2y_{-1})\epsilon+\dots\label{eq:B25}.
\end{align}
We now apply (\ref{eq:B22})-(\ref{eq:B25}) to (\ref{eq:40}) and multiply by $\epsilon,$ yielding
\begin{align}
    0 &= \epsilon^3\left(y_{-1}^2\epsilon^{-2}+(2y_{-1}y_0+2x_0y_{-1})\epsilon^{-1} \right. \nonumber\\
    &\qquad+\left(x_0^2 +x_0y_0+x_1y_{-1}+y_0^2+2y_{-1}y_1\right) \nonumber\\
    &\qquad \left. +2\left(x_0x_1+x_0y_1+x_1y_0+x_2y_{-1}+y_{-1}y_2+y_0y_1\right) \epsilon+\dots \right) \nonumber\\
    &- \epsilon^2\left(\mu y_{-1}^2\epsilon^{-2}+ (\lambda y_{-1}+\mu(x_0y_{-1} +2y_{-1}y_0))\epsilon^{-1}\right. \nonumber\\
    &\qquad+\lambda y_0 +\mu(x_0y_0+x_1y_{-1}+y_0^2+2y_{-1}y_1) \nonumber \\
    &\qquad  + (\lambda y_1+\mu\left(x_0y_1+x_1y_0+x_2y_{-1}+2y_{-1}y_2+2y_0y_1\right))\epsilon \nonumber\\
    &\qquad\left.+ (\lambda y_2+\mu\left(x_0y_2+x_1y_1+x_2y_0+x_3y_{-1}+y_1^2+2y_{-1}y_3+2y_0y_2\right))\epsilon^2\right. \nonumber\\
    &\qquad\left. + \dots\right)  \nonumber\\
    &+ \epsilon\lambda\sigma \left(x_0y_{-1}\epsilon^{-1} + (x_0y_0+x_1y_{-1}) +(x_0y_1+x_1y_0+x_2y_{-1})\epsilon \right.\nonumber\\
    &\qquad+ (x_0y_2+x_1y_1+x_2y_0+x_3y_{-1})\epsilon^2  \nonumber\\
    &\qquad +\left. (x_0y_3+x_1y_2+x_2y_1+x_3y_0+x_4y_{-1})\epsilon^3 \dots\right). \label{eq:B26}
\end{align}
Equating the $\mathcal{O}(1)$ terms to zero, we have
\begin{equation}
    0 = \lambda\sigma x_0y_{-1}-\mu y_{-1}^2 = x_0^2(\lambda\sigma \psi_{-1}-\mu \psi_{-1}^2), \label{eq:B27}
\end{equation}
and thus $x_0 = y_{-1}=0.$ Equating the $\mathcal{O}(\epsilon)$ terms to zero, we have
\begin{equation}
    0 = y_{-1}^2+\lambda y_{-1}+\mu(x_0y_{-1}+2y_{-1}y_0)+\lambda\sigma(x_0y_0+x_1y_{-1}),\label{eq:B28}
\end{equation}
which is seen to be trivially satisfied as a result of (\ref{eq:B27}). Therefore, we look to the $\mathcal{O}(\epsilon^2)$ terms to determine $x_1.$ Equating those coefficients to zero leads to 
\begin{align}
    0 &= 2y_{-1}y_0+2x_0y_{-1}-\lambda y_0 -\mu(x_0y_0+x_1y_{-1}+y_0^2+2y_{-1}y_1)\nonumber\\
    &\qquad +\lambda\sigma(x_0y_1+x_1y_0+x_2y_{-1}) \nonumber \\
    &= -\lambda y_0 - \mu y_0^2+\lambda\sigma x_1y_0\nonumber\\
    &= -\psi_{-1}(x_1-\phi_1)(\lambda-\mu\psi_{-1}\phi_1+(\mu\psi_{-1}-\lambda\sigma)x_1). \label{eq:B29}
\end{align}
Of the two solutions to this equation, we are interested in $x_1 = \phi_1 = 1/\sigma,$ which in turn implies that $y_1 = 0$.

Looking now for $x_2,$ we equate the $\mathcal{O}(\epsilon^3)$ coefficients to zero:
\begin{align}
    0 &= x_0^2+x_0y_0+x_1y_{-1}+y_0^2+2y_{-1}y_1 \nonumber\\
    &\qquad-\lambda y_1 -\mu(x_0y_1+x_1y_0+x_2y_{-1}+2y_{-1}y_2+2y_0y_1)\nonumber\\
    &\qquad+(x_0y_2+x_1y_1+x_2y_0+x_3y_{-1}).\label{eq:B30}
\end{align}
which is also trivially satisfied as all terms either cancel with another or contain a factor of $x_0,y_{-1},$ or $y_0.$ Thus, we turn to $\mathcal{O}(\epsilon^4)$ to determine $x_2.$ Equating the coefficients to zero gives
\begin{align}
    0 &= 2(x_0x_1+x_0y_1+x_1y_0+x_2y_{-1}+y_{-1}y_2+y_0y_1)-\lambda y_2 \nonumber\\
    &\qquad -\mu\left(x_0y_2+x_1y_1+x_2y_0+x_3y_{-1}+y_1^2+2y_{-1}y_3+2y_0y_2\right) \nonumber\\
    &\qquad +\lambda\sigma(x_0y_3+x_1y_2+x_2y_1+x_3y_0+x_4y_{-1}) \nonumber \\
    &= y_1(-\mu(x_1+y_1)+\lambda\sigma x_2). \label{eq:B31}
\end{align}
The solution we're interested in for $x_2$ comes from $y_1=0,$ which can be expressed in terms of $x_2$ as 
\begin{equation}
    0 = \psi_{-1}(x_2-\phi_2), \label{eq:B32}\\
\end{equation}
and thus
\begin{equation}
    x_2  =\phi_2 = \frac{\delta_c+\mu-\sigma}{\lambda\sigma^2}.\label{eq:B33}
\end{equation}
At this point, we have a second order expansion of the approximate equilibrium $x^*:$\
\begin{equation}
x^* \approx \frac{1}{\sigma}\epsilon + \frac{\delta_c+\mu-\sigma}{\lambda\sigma^2}\epsilon^2+\mathcal{O}(\epsilon^3).\label{eq:B34}
\end{equation}
Now with the relation $w^* = \frac{\sigma}{\epsilon}x^*,$ we conclude that
\begin{equation}
    w^* \approx 1+\frac{\delta_c+\mu-\sigma}{\lambda\sigma}\epsilon+\mathcal{O}(\epsilon^2).\label{eq:B35}
\end{equation}
\end{appendices}

\bibliographystyle{spbasic}      % basic style, author-year citations
\bibliography{references.bib}   % name your BibTeX data base

\end{document}